\documentclass[preprint,11pt]{elsarticle}

\newtheorem{theorem}{Theorem}[section]
\newtheorem{lemma}{Lemma}[section]
\newtheorem{corollary}{Corollary}[section]
\newtheorem{proposition}{Proposition}[section]

\newtheorem{remark}{Remark}[section]
\newcommand{\ignore}[1]{}{}

\usepackage{amssymb}
\usepackage{amsmath}

\usepackage{color}
\usepackage{soul}
\usepackage[dvipsnames]{xcolor}

\usepackage[colorlinks=true, linkcolor=blue, citecolor=blue]{hyperref}

\def\1{{{\mbox{${\rm{1\negthinspace\negthinspace I}}$}}}}

\newcommand\beq{\begin{equation}}
\newcommand\eeq{\end{equation}}

\usepackage{ifthen}
\usepackage{xkeyval}
\usepackage{todonotes}
\begin{document}

\begin{frontmatter}

\title{Self-normalized Cram\'{e}r  moderate deviations for a supercritical Galton-Watson process}
\author[cor1]{Xiequan Fan }
\author[cor3]{Qi-Man Shao }
\address[cor1]{Center for Applied Mathematics,
Tianjin University, Tianjin 300072,  China}
\address[cor3]{Department of Statistics and Data Science, Southern University of Science and Technology,
 Shenzhen  518000, China.}


\begin{abstract}
Let $(Z_n)_{n\geq0}$ be a supercritical Galton-Watson process.
Consider the   Lotka-Nagaev estimator  for the offspring mean.
In this paper, we establish self-normalized Cram\'{e}r type moderate deviations and Berry-Esseen's bounds  for  the Lotka-Nagaev estimator.
 The results are believed to be optimal or near optimal.
\end{abstract}

\begin{keyword}   Lotka-Nagaev estimator; offspring mean;  Self-normalized processes;
 Cram\'{e}r moderate deviations; Berry-Esseen's bounds
\vspace{0.3cm}
\MSC primary 60J80; 60F10;   secondary  62F03; 62F12
\end{keyword}

\end{frontmatter}




\section{Introduction}
\setcounter{equation}{0}
A Galton-Watson process   can be described as follows
 \begin{equation} \label{GWP}
 Z_0=1,\ \ \ \ Z_{n+1}= \sum_{i=1}^{Z_n} X_{n,i}, \   \ \ \ \textrm{for } n \geq 0,
 \end{equation}
where $X_{n,i}$ is the offspring number of   the $i$-th individual of the generation $n.$
Moreover, the random variables $  (X_{n,i})_{i\geq 1} $ are independent of each other with common  distribution law
  \begin{equation} \label{ssfs}
 \mathbb{P}(X_{n,i} =k   )     = p_k,\ \ \ k \in \mathbb{N},
\end{equation}
  and are also independent of $Z_n.$

An important task in statistical inference of Galton-Watson processes is to estimate  the average offspring number  of an individual $m,$    usually termed the offspring mean. Clearly, it holds  $$m=\mathbb{E}Z_1=\mathbb{E} X_{n,i} =\sum_{k=0}^\infty  k   p_k.$$
 Denote  $v$   the standard variance of $Z_1$, that is
\begin{eqnarray}\label{defv}
  \upsilon^2=\mathbb{E} (Z_1-m)^2.
\end{eqnarray}
To avoid triviality,  assume that  $v> 0.$
 For estimation of the offspring mean $m$,  the Lotka-Nagaev    \cite{L39,N67}   estimator $Z_{n+1}/Z_{n}$  plays an important role. For the Galton-Watson processes, Athreya \cite{A94} has established large deviations
for the  normalized Lotka-Nagaev estimator (see also Chu \cite{C18}  for self-normalized large deviations);
Ney and Vidyashankar \cite{NV03,NV04} obtained sharp rate estimates for the large deviation behavior of  the Lotka-Nagaev  estimator;
Bercu and Touati \cite{BT08} proved an exponential inequalities for the Lotka-Nagaev estimator via self-normalized martingale
method.  
The main purpose of this paper is to establish self-normalized Cram\'{e}r moderate deviations  for the Lotka-Nagaev   estimator $Z_{n+1}/Z_{n}$  for the Galton-Watson processes.

 The paper is organized as follows. In Section \ref{sec2}, we present
 Cram\'{e}r moderate deviations  for the self-normalized  Lotka-Nagaev  estimator,
provided that $(Z_n)_{n\geq0}$ or $(X_{n,i})_{1\leq i \leq Z_n}$ can be observed.
  In Section  \ref{sec3}, we  present some applications of our results in statistics.
 The rest sections devote to the proofs of  theorems.


\section{Main results}\label{sec2}
\setcounter{equation}{0}

\subsection{ $(Z_{k})_{k\geq 0}$ can be observed}
Assume that the total populations $(Z_{k})_{k\geq 0}$ of all generations can be observed.
For any $n_0\geq 0,$ we define
\begin{eqnarray}
M_{n_0,n}= \frac{ \sum_{k=n_0}^{n_0+n-1}    \sqrt{Z_{ k}} (  \frac{Z_{ k+1}}{Z_{ k}} -m )}{\sqrt{\sum_{k=n_0}^{n_0+n-1}   Z_{  k}  (  \frac{Z_{ k+1}}{Z_{ k}}  - m  )^2 } \ }.
\end{eqnarray}
 We  assume that the set of extinction of the process
$(Z_{k})_{k\geq 0}$ is negligible with respect to the annealed law $\mathbb{P}$. Then $M_{n_0,n}$ is well defined $\mathbb{P}$-a.s.
 As $(Z_k)_{k=n_0,...,n_0+n}$ can be observed, $M_{n_0,n}$ can be regarded as  a time type self-normalized process for the Lotka-Nagaev estimator $Z_{k+1}/Z_{k}$.
The following theorem gives a self-normalized Cram\'{e}r moderate deviation result for the Galton-Watson processes.

\begin{theorem}\label{th01}
Assume that     $ \mathbb{E} Z_1 ^{2+\rho}< \infty$ for some  $  \rho \in (0,  1]$.
\begin{description}
  \item[\textbf{[i]}]  If $\rho \in (0, 1)$, then for all $x \in [0,\   o(  \sqrt{n} )),$
\begin{equation}\label{dfdsf12}
\bigg|\ln\frac{\mathbb{P}(M_{n_0,n} \geq x)}{1-\Phi(x)} \bigg| \leq  C_\rho \bigg(    \frac{ x^{2+\rho} }{n^{\rho/2}}  +   \frac{ (1+x)^{1-\rho(2+\rho)/4} }{n^{\rho(2-\rho)/8}}     \bigg) ,
\end{equation}
where $C_\rho$ depends only on  the constants $\rho,  v$ and $ \mathbb{E} Z_1 ^{2+\rho}$.
  \item[\textbf{[ii]}]  If $\rho =1$, then for all $x \in [0,\   o(  \sqrt{n} )),$
\begin{equation}\label{dfdsf13}
\bigg|\ln\frac{\mathbb{P}(M_{n_0,n} \geq x)}{1-\Phi(x)} \bigg| \leq  C \bigg(    \frac{ x^{3} }{\sqrt{n} }  + \frac{ \ln n }{\sqrt{n} }  +  \frac{ (1+x)^{1/4} }{n^{1/8}}     \bigg) ,
\end{equation}
where $C$ depends only on  the constants $  v$ and $ \mathbb{E} Z_1 ^{3}$.
\end{description}
In particular, the   inequalities \eqref{dfdsf12} and \eqref{dfdsf13} together implies  that
\begin{eqnarray}\label{dsfaf}
\frac{\mathbb{P}(M_{n_0,n} \geq x)}{1-\Phi(x)} =1+o(1)
\end{eqnarray}
uniformly for $n_0 \in \mathbb{N}$ and for $ x \in [0, \, o(n^{\rho/(4+2\rho)}))$ as $n\rightarrow \infty$.
Moreover, the same inequalities remain valid when $\frac{\mathbb{P}( M_{n_0,n} \geq x)}{1-\Phi \left( x\right)}$ is replaced  by $\frac{\mathbb{P}(M_{n_0,n} \leq -x)}{ \Phi \left( -x\right)}$.
\end{theorem}

Notice that $C_\rho$ and $C$ do not depend on $n_0.$ Thus  \eqref{dsfaf}  holds uniformity in $n_0$, which is of particular interesting in applications. For instance,  due to the uniformity, in \eqref{dsfaf}
we can take $n_0$  as a function of $n$.

Equality (\ref{dsfaf}) implies that $\mathbb{P}(M_{n_0,n} \leq x) \rightarrow  \Phi(x)$ as  $n$ tends to $\infty$. Thus Theorem \ref{th01} implies the central limit theory
for $M_{n_0,n}$. Moreover,  equality (\ref{dsfaf}) states that the relative error of normal approximation for $M_{n_0,n}$ tends to zero   uniformly for $  x \in [0, \, o(n^{\rho/(4+2\rho)})) $ as $n\rightarrow \infty$.

Theorem  \ref{th01}   implies the following moderate deviation principle  (MDP) result for the  time type self-normalized Lotka-Nagaev estimator.
\begin{corollary} \label{coro01}
Assume the conditions  of Theorem  \ref{th01}.
Let $(a_n)_{n\geq1}$ be any sequence of real numbers satisfying $a_n \rightarrow \infty$ and $a_n/ \sqrt{n}    \rightarrow 0$
as $n\rightarrow \infty$.  Then  for each Borel set $B$,
\begin{eqnarray}
- \inf_{x \in B^o}\frac{x^2}{2}  \leq   \liminf_{n\rightarrow \infty}\frac{1}{a_n^2}\ln \mathbb{P}\bigg(  \frac{M_{n_0,n} }{ a_n  }     \in B \bigg)
  \leq \limsup_{n\rightarrow \infty}\frac{1}{a_n^2}\ln \mathbb{P}\bigg(\frac{ M_{n_0,n} }{ a_n }  \in B \bigg) \leq  - \inf_{x \in \overline{B}}\frac{x^2}{2}   ,   \label{SMDP}
\end{eqnarray}
where $B^o$ and $\overline{B}$ denote the interior and the closure of $B$, respectively.
\end{corollary}

\begin{remark}
From (\ref{dfdsf12}) and \eqref{dfdsf13},   it is easy to derive the following Berry-Esseen bound  for the self-normalized Lotka-Nagaev estimator:
\begin{eqnarray}
\Big|\mathbb{P}(M_{n_0,n} \leq x) -    \Phi(x) \Big| \leq  \frac{   C_\rho }{n^{\rho(2-\rho)/8}},
\end{eqnarray}
where $C_\rho$ depends only on  the constants $\rho,  v$ and $ \mathbb{E} Z_1 ^{2+\rho}$. When  $\rho> 1$,
by the self-normalized Berry-Esseen bound  for martingales in Fan and Shao \cite{F18}, we can get
 a Berry-Esseen bound of order $n^{-\frac{\rho}{6+2\rho}}$.
\end{remark}

The last remark gives a self-normalized  Berry-Esseen bound for the Lotka-Nagaev estimator, while the next theorem presents
a  normalized  Berry-Esseen bound for the Lotka-Nagaev estimator.
Denote
$$H_{n_0, n}= \frac{1}{ \sqrt{n } v }\sum_{k=n_0}^{n_0+n-1}    \sqrt{Z_{ k}} \Big ( \frac{Z_{k+1}}{Z_{k}} -m \Big). $$
Notice that the random variables $(X_{k,i})_{1\leq i\leq Z_k}$
have the   same distribution as $Z_1,$ and that $(X_{k,i})_{1\leq i\leq Z_k}$ are independent of $Z_k$. Then for the Galton-Watson processes, it holds
$$\mathbb{E} [ (Z_{k+1}  -m Z_{k})^2 | Z_k] = \mathbb{E} [ ( \sum_{i=1}^{Z_k}  (X_{k,i} -m) )^2 | Z_k] = Z_{k}  \upsilon^2.$$
It is easy to see that $H_{n_0, n}=  \sum_{k=n_0}^{n_0+n-1}  \frac{1}{ \sqrt{\, n v^2 / Z_{ k}  }}   \Big ( \frac{Z_{k+1}}{Z_{k}} -m \Big)  .$
Thus $H_{n_0,n}$  can be regarded as a normalized process for the Lotka-Nagaev estimator $Z_{k+1}/Z_{k}$.  We  have the following normalized Berry-Esseen bounds
for   the Galton-Watson processes.
\begin{theorem} \label{jsdf}
Assume the conditions of Theorem  \ref{th01} are satisfied.
\begin{description}
  \item[\textbf{[i]}]  If $\rho \in (0, 1)$, then
 \begin{eqnarray}\label{INES29}
\sup_{x \in \mathbb{R}}\Big|\mathbb{P}( H_{n_0,n} \leq x) -  \Phi(x)  \Big| \leq \frac{ C_\rho }{ n^{\rho/2}},
\end{eqnarray}
where $C_\rho$ depends only on $\rho, v$ and $ \mathbb{E} Z_1 ^{2+\rho}$.
  \item[\textbf{[ii]}]  If $\rho =1$, then
 \begin{eqnarray}\label{INES210}
\sup_{x \in \mathbb{R}}\Big|\mathbb{P}( H_{n_0,n} \leq x) -  \Phi(x)  \Big| \leq C   \frac{ \ln n }{ \sqrt{n} },
\end{eqnarray}
where $C$ depends only  on $ v$ and $ \mathbb{E} Z_1 ^{3}$.
 \end{description}
 Moreover, the same inequalities remain  valid when $H_{n_0,n}$ is replaced  by $-H_{n_0,n}$.
\end{theorem}

The convergence rates of (\ref{INES29}) and (\ref{INES210}) are same to the best possible convergence rates of the Berry-Esseen bounds for martingales,
 see  Theorem 2.1   of  Fan \cite{F19} and its comment.
Notice that $H_{n_0,n}$ is a martingale with respect  to the natural filtration.


\subsection{ $(X_{n,i})_{1\leq i \leq Z_n}$ can be observed for some $n$}

Assume that the offspring numbers $(X_{n,i})_{1\leq i \leq Z_n}$ of each individual in some generation  $n$   can be observed.
 Denote
$$T_n= \frac{ Z_n \, ( \frac{Z_{n+1}}{Z_{n}} -m  )}{\sqrt{\sum_{i=1}^{Z_{n}}(X_{n,i}- \frac{Z_{n+1}}{Z_{n}})^2}} $$
the space type self-normalized process  for the Lotka-Nagaev estimator $Z_{n+1}/Z_{n}$. The following theorem gives a Cram\'{e}r moderate deviation result for the  space type self-normalized  Lotka-Nagaev estimator $T_n$.
\begin{theorem}\label{fsfs}
Assume that  $p_0=0$ and   $ \mathbb{E} Z_1 ^{2+\rho}< \infty$ for some  $  \rho \in (0,  1]$. Then
\begin{eqnarray}\label{fsdcv01}
\bigg|\ln\frac{\mathbb{P} (T_{n} \geq x)}{1-\Phi(x)} \bigg| =  O\Big( \frac{1+x^{2+\rho}}{ n^{\rho/2}}\Big)
\end{eqnarray}
uniformly for $x \in [0, \, o(\sqrt{n}) )$ as $n\rightarrow \infty$.
Moreover, the same  equality remains valid when $\frac{\mathbb{P} (T_n \geq x)}{1-\Phi \left( x\right)}$ is replaced  by $\frac{\mathbb{P}  (T_n \leq -x)}{ \Phi \left( -x\right)}$.
\end{theorem}

The condition $p_0=0$ means that   each individual has at least one offspring. Moreover,  it also implies that
$Z_n \rightarrow \infty$ a.s.\ as $n\rightarrow \infty.$ Then by law of large numbers, we have $\frac{Z_{n+1}}{Z_{n}}$ tends
to $m$ a.s.\ as $n\rightarrow \infty.$

For the Galton-Watson processes,  we refer to \cite{A94} for closely related results of Theorem \ref{fsfs}, where Athreya
has established a precise large deviation rate for the Lotka-Nagaev estimator $Z_{n+1}/Z_{n}$.

Using the inequality $|e^x-1| \leq e^C |x|$ valid for $|x| \leq C,$ from  Theorem \ref{fsfs}, we obtain the following
estimation for the relative error of normal approximation.
\begin{corollary}\label{corollary01}
Assume the conditions  of Theorem \ref{fsfs}. Then
\begin{equation}\label{dfrezero}
\frac{\mathbb{P} \big( T_n \geq x  \big)}{1-\Phi(x)} =1+ O\Big(\frac{1+x^{2+\rho}}{n^{\rho/2}}\Big)
\end{equation}
uniformly for   $x \in [0, \, O( n^{\rho/(4+2\rho)} ))$ as $n\rightarrow \infty$.
In particular, it implies that
\begin{equation}\label{rezero}
\frac{\mathbb{P} \big( T_n \geq x  \big)}{1-\Phi(x)} =1+ o(1)
\end{equation}
uniformly for $x \in [0, \, o( n^{\rho/(4+2\rho)} ))$ as $n\rightarrow \infty$.
Moreover, the same  equalities remain  valid when $T_n $ is replaced  by $-T_n$.
\end{corollary}

Inequality (\ref{rezero}) implies  that the relative error of normal approximation for $T_n$ tends to zero   uniformly for $x \in [0, \, o(  n^{\rho/(4+2\rho)} )).$
Clearly, the range of validity  for (\ref{rezero})  coincides  with the self-normalized Cram\'{e}r moderate deviation result of Shao \cite{S99}  for iid random variables.

By an argument similar to the proof of Corollary  \ref{coro01},
Theorem  \ref{fsfs} also  implies  the following self-normalized   MDP  result.
\begin{corollary}\label{corollary02}
Assume the conditions   of Theorem \ref{fsfs}.
Let $(a_n)_{n\geq1}$ be any sequence of real numbers satisfying $a_n \rightarrow \infty$ and $a_n/ \sqrt{n}    \rightarrow 0$
as $n\rightarrow \infty$.  Then  for each Borel set $B$,
\begin{eqnarray}
- \inf_{x \in B^o}\frac{x^2}{2}  \leq   \liminf_{n\rightarrow \infty}\frac{1}{a_n^2}\ln \mathbb{P} \bigg(  \frac{T_n }{ a_n  }     \in B \bigg)
  \leq \limsup_{n\rightarrow \infty}\frac{1}{a_n^2}\ln \mathbb{P} \bigg(\frac{ T_n }{ a_n }  \in B \bigg) \leq  - \inf_{x \in \overline{B}}\frac{x^2}{2}   ,   \label{MDP}
\end{eqnarray}
where $B^o$ and $\overline{B}$ denote the interior and the closure of $B$, respectively.
\end{corollary}

From  Theorem \ref{fsfs}, we get the following self-normalized Berry-Esseen bound for $T_n$.
\begin{corollary}\label{corollary03} Assume the conditions   of Theorem \ref{fsfs}.  Then
\begin{eqnarray}\label{ths01}
\sup_{x \in \mathbb{R}}\Big|\mathbb{P} ( T_{n} \leq x) -  \Phi(x)  \Big| \leq \frac{ C_\rho }{ n^{\rho/2}},
\end{eqnarray}
where $C_\rho$ does not depend on $n.$
\end{corollary}

Clearly, the convergence rate for the Berry-Esseen bound of Corollary \ref{corollary03} is consistent with the classical case of iid random variables (cf.\ Bentkus and  G\"{o}tze \cite{BG96}), and therefore it is optimal under the stated conditions.

\begin{remark}
Following  the proof of Theorem  \ref{fsfs}, the results (\ref{fsdcv01})-(\ref{ths01}) remain true when $T_n$ is replaced by
$$  \widetilde{T}_n= \frac{ Z_n \, ( Z_{n+1}/Z_{n} -m  )}{\sqrt{\sum_{i=1}^{Z_{n}}(X_{n,i}-m  )^2}}. $$
\end{remark}

\section{Applications}\label{sec3}
\setcounter{equation}{0}
Cram\'{e}r  moderate deviations certainly have a lot of  applications in statistics.
 \subsection{$p$-value for hypothesis testing}
 Self-normalized Cram\'{e}r  moderate deviations can be applied to hypothesis testing  of $m$  for the Galton-Watson processes.
 When $(Z_{k})_{k=n_0,...,n_0+n}$ can be observed, we can make use of  Theorem \ref{th01}
to estimate $p$-value.
Assume that $ \mathbb{E} Z_1 ^{2+\rho}< \infty$ for some $0 < \rho \leq 1$, and  that $m> 1$.  Let $(z_{k})_{k=n_0,...,n_0+n}$ be an observation of $(Z_{k})_{k=n_0,...,n_0+n}$.
In order to estimate the offspring mean $m$,
we can make use of the Harris estimator  \cite{BT08} given by
$$\widehat{m}_n=\frac{\sum_{k=n_0}^{n_0+n-1}  Z_{ k+1}}{\sum_{k=n_0}^{n_0+n-1}Z_{ k}}. $$
Then observation for the Harris estimator is
$$\widehat{m}_n=\frac{\sum_{k=n_0}^{n_0+n-1}  z_{ k+1}}{\sum_{k=n_0}^{n_0+n-1}z_{ k}}. $$
By Theorem \ref{th01}, it is easy to see that
\begin{equation}   \label{sfdfs}
\frac{\mathbb{P}( M_{n_0,n}  \geq x)}{1-\Phi \left( x\right)}=1+o(1)\ \ \ \ \textrm{and} \ \ \ \ \frac{\mathbb{P}(M_{n_0,n}  \leq- x)}{1-\Phi \left( x\right)}=1+o(1)
\end{equation}
uniformly for $ x \in [0,  o (  n^{\rho/(4+2\rho)}  )).$   Notice that $ 1-\Phi \left( x\right) =  \Phi \left( -x\right). $ Thus, by (\ref{sfdfs}), the probability  $\mathbb{P}(M_{n_0,n}   >  |\widetilde{m}_n|)$  is almost equal to $2 \Phi \left( -|\widetilde{m}_n|\right) $,   where
$$\widetilde{m}_n= \frac{ \sum_{k=n_0}^{n_0+n-1}    \sqrt{z_{ k}} ( z_{ k+1}/z_{ k} -\widehat{m}_n )}{\sqrt{\sum_{k=n_0}^{n_0+n-1}   z_{  k}  ( z_{ k+1}/z_{ k} - \widehat{m}_n  )^2 } \ }.  $$

\subsection{Construction of confidence intervals }
\subsubsection{The data $(Z_{k})_{k\geq 0}$ can be observed}
Cram\'{e}r  moderate deviations  can  be also applied to  construction of confidence intervals of $m$. We make use of Theorem \ref{th01}
to construct confidence intervals.
\begin{proposition}\label{c0klfdgs}
Assume that $ \mathbb{E} Z_1 ^{2+\rho}< \infty$ for  some $  \rho \in (0,  1]$. Let $\kappa_n \in (0, 1).$  Assume that
\begin{eqnarray}\label{keldffdet}
 \big| \ln \kappa_n \big| =o \big( n^{\rho/(2+\rho)}  \big) .
\end{eqnarray}
Let
\begin{eqnarray*}
a_{n_0, n}&=& \Big(  \sum_{k=n_0}^{n_0+n-1} \sqrt{Z_{k}} \Big)^2- \big(\Phi^{-1}(1-\kappa_n /2)\, \big)^2 \sum_{k=n_0}^{n_0+n-1} Z_{k}, \\
b_{n_0, n}&=& 2 \big(\Phi^{-1}(1-\kappa_n /2) \big)^2 \sum_{k=n_0}^{n_0+n-1} Z_{k+1} - 2 \Big( \sum_{k=n_0}^{n_0+n-1} \frac{Z_{k+1}}{\sqrt{Z_{k }}} \Big)\Big( \sum_{k=n_0}^{n_0+n-1} \sqrt{Z_{k}}  \Big), \\
c_{n_0, n}&=& \Big(  \sum_{k=n_0}^{n_0+n-1} \frac{Z_{k+1}}{\sqrt{Z_{k }}}  \Big)^2 - \big(\Phi^{-1}(1-\kappa_n /2) \big)^2 \sum_{k=n_0}^{n_0+n-1} \frac{Z_{k+1}^2}{Z_{k}}.
\end{eqnarray*}
Then $[A_{n_0, n},B_{n_0, n}]$, with
\begin{eqnarray*}
A_{n_0, n}=\frac{- b_{n_0, n} - \sqrt{b_{n_0, n}^2-4 a_{n_0, n}c_{n_0, n}  }}{ 2 a_{n_0, n} }
\end{eqnarray*}
and
\begin{eqnarray*}
B_{n_0, n}=\frac{ -b_{n_0, n} + \sqrt{b_{n_0, n}^2-4 a_{n_0, n}c_{n_0, n}  }}{ 2 a_{n_0, n} },
\end{eqnarray*}
is a  $1-\kappa_n$ confidence interval for $m$, for $n$ large enough.
%
\end{proposition}
\emph{Proof.}  Notice that $ 1-\Phi \left( x\right) =  \Phi \left( -x\right). $ Theorem \ref{th01} implies that
\begin{equation}   \label{sfdfsdfs}
\frac{\mathbb{P}( M_{n_0,n}  \geq x)}{1-\Phi \left( x\right)}=1+o(1)\ \ \ \ \textrm{and} \ \ \ \ \frac{\mathbb{P}(M_{n_0,n}  \leq- x)}{1-\Phi \left( x\right)}=1+o(1)
\end{equation}
uniformly for $0\leq x=o (   n^ {\rho/(4+2\rho)}   )$, see \eqref{dsfaf}.  When $\kappa_n$ satisfies the condition (\ref{keldffdet}),
  the upper $(\kappa_n/2)$th quantile of a standard normal distribution satisfies
$$\Phi^{-1}( 1-\kappa_n/2) = O(\sqrt{| \ln \kappa_n |} ), $$  which  is of order $o\big( n^ {\rho/(4+2\rho)}  \big).$
Then applying (\ref{sfdfsdfs}) to the last equality, we complete the proof of  Proposition \ref{c0klfdgs}.
Notice that $A_{n_0, n}$ and $B_{n_0, n}$ are solutions of the following equation
$$  \frac{ \sum_{k=n_0}^{n_0+n-1}    \sqrt{Z_{ k}} ( Z_{ k+1}/Z_{ k} -x )}{\sqrt{\sum_{k=n_0}^{n_0+n-1}   Z_{  k}  ( Z_{ k+1}/Z_{ k} - x  )^2 } \ }=\Phi^{-1}(1-\kappa_n /2).  $$
This completes the proof of Proposition \ref{c0klfdgs}. \hfill\qed

\subsubsection{The data $(X_{n,i})_{1\leq i \leq Z_n}$ can be observed}
When $(X_{n,i})_{1\leq i \leq Z_n}$ can be observed, we can make use of  Corollary \ref{corollary01}
to construct confidence intervals.
\begin{proposition}\label{c0kls}
Assume that $ \mathbb{E} Z_1 ^{2+\rho}< \infty$ for some $  \rho \in (0,  1]$. Let $\kappa_n \in (0, 1).$  Assume that
\begin{eqnarray}\label{keldet}
 \big| \ln \kappa_n \big| =o \big( n^{\rho/(2+\rho)}  \big) .
\end{eqnarray}
Let $$\Delta_n=\frac{\Phi^{-1}(1-\kappa_n/2)  }{Z_n } \sqrt{\sum_{i=1}^{Z_{n}}(X_{n,i}-\frac{Z_{n+1}}{Z_{n}})^2}.$$
Then   $[A_n,B_n]$, with
\begin{eqnarray*}
A_n=\frac{ Z_{n+1} }{ Z_{n}}  -\Delta_n  \quad   \textrm{and} \ \quad  
B_n=\frac{ Z_{n+1} }{ Z_{n}} +\Delta_n, 
\end{eqnarray*}
is a  $1-\kappa_n$ confidence interval for $m$, for $n$ large enough.
\end{proposition}
\emph{Proof.} Corollary \ref{corollary01} implies that
\begin{equation} \label{tphisnds4}
\frac{\mathbb{P} (T_n  \geq x)}{1-\Phi \left( x\right)}=1+o(1)\ \ \ \ \textrm{and} \ \ \ \ \frac{\mathbb{P} (T_n  \leq- x)}{1-\Phi \left( x\right)}=1+o(1)
\end{equation}
uniformly for $0\leq x=o (   n^ {\rho/(4+2\rho)}   ).$  When $\kappa_n$ satisfies the condition (\ref{keldffdet}),
  the upper $(\kappa_n/2)$th quantile of a standard normal distribution satisfies
 $\Phi^{-1}( 1-\kappa_n/2) = O(\sqrt{| \ln \kappa_n |} ), $  which  is of order $o\big( n^ {\rho/(4+2\rho)}  \big).$
Then applying (\ref{tphisnds4}) to the last equality, we complete the proof of  Proposition \ref{c0kls}. \hfill\qed

When the risk probability $\kappa_n$ goes to $0$, we have the following more general result.
\begin{proposition}\label{zfgsgg} Assume that $ \mathbb{E} Z_1 ^{2+\rho}< \infty$ for some $  \rho \in (0,  1]$. Let $\kappa_n \in (0, 1)$ such that $k_n \rightarrow 0$. Assume that
\begin{eqnarray}
 \big| \ln \kappa_n \big| =o \big(  \sqrt{n}   \big) .
\end{eqnarray}
Let $$\Delta_n=\frac{\sqrt{ 2 |\ln (\kappa_n/2)|}  }{Z_n } \sqrt{\sum_{i=1}^{Z_{n}}(X_{n,i}-\frac{Z_{n+1}}{Z_{n}})^2}.$$
Then   $[A_n,B_n]$,   with
\begin{eqnarray*}
A_n=\frac{ Z_{n+1} }{ Z_{n}}  -\Delta_n  \quad   \textrm{and} \ \quad  
B_n=\frac{ Z_{n+1} }{ Z_{n}} +\Delta_n, 
\end{eqnarray*}
is a  $1-\kappa_n$ confidence interval for $m $, for $n$ large enough.
\end{proposition}
\emph{Proof.} By Theorem \ref{fsfs}, we have
\begin{equation}\label{ggsg1d2}
\frac{\mathbb{P} (T_n  \geq x)}{1-\Phi \left( x\right)}=\exp\Big\{ \theta C  \frac{1+x^{2+\rho}}{ n^{\rho/2}} \Big \} \ \ \ \  \textrm{and} \ \ \ \ \frac{\mathbb{P} (T_n \leq - x)}{1-\Phi \left( x\right)}=\exp\Big\{ \theta C  \frac{1+x^{2+\rho}}{ n^{\rho/2}} \Big \}
\end{equation}
uniformly for $0\leq x=o (   \sqrt{n}   ),$ where $\theta \in [-1, 1]$. Notice that $$1-\Phi \left( x_n\right) \sim \frac{1}{x_n\sqrt{2\pi}}e^{-x_n^2/2}= \exp\bigg\{-\frac{x_n^2}{2}\Big(1+\frac{2}{x_n^2}\ln (x_n\sqrt{2\pi})   \Big) \bigg\} ,\ x_n \rightarrow \infty. $$
  Since $k_n \rightarrow 0$, the last line implies that
the upper $(\kappa_n/2)$th quantile of the distribution $$1-\Big(1-\Phi \left( x\right)\Big)\exp\Big\{ \theta C  \frac{1+x^{2+\rho}}{ n^{\rho/2}} \Big\}$$
converges to  $ \sqrt{ 2 |\ln (\kappa_n/2)|}$, which  is of order $o\big(\sqrt{n} \big)$  as $n\rightarrow \infty.$
Then applying (\ref{ggsg1d2}) to $T_n$, we complete the proof of  Proposition \ref{zfgsgg}. \hfill\qed

\subsubsection{The parameter $\upsilon^2$ is known}
When $\upsilon^2$  is known,  we can apply normalized Berry-Esseen bounds (cf.\ Theorem \ref{jsdf})
to construct confidence intervals.
\begin{proposition} \label{pro3.3}
Assume that $ \mathbb{E} Z_1 ^{2+\rho}< \infty$ for some $  \rho \in (0,  1]$.    Let $\kappa_n \in (0, 1).$  Assume that
\begin{eqnarray}  \label{fs25sv}
 \big| \ln \kappa_n \big| =o \big(  \log n   \big) .
\end{eqnarray}
Then   $[A_n,B_n]$,   with
\begin{eqnarray*}
A_n=\frac{ \sum_{k=n_0}^{n_0+n } Z_{k+1}/\sqrt{Z_k}-\sqrt{n } v \Phi^{-1}(1-\kappa_n/2)   }{ \sum_{k=n_0}^{n_0+n } \sqrt{Z_{k}}}
\end{eqnarray*}
and
\begin{eqnarray*}
B_n=\frac{ \sum_{k=n_0}^{n_0+n } Z_{k+1}/\sqrt{Z_k}+\sqrt{n } v \Phi^{-1}(1-\kappa_n/2)  }{ \sum_{k=n_0}^{n_0+n } \sqrt{Z_{k}}} ,
\end{eqnarray*}
is a  $1-\kappa_n$ confidence interval for $m$, for $n$ large enough.
\end{proposition}
\emph{Proof.} Theorem \ref{jsdf}   implies that
\begin{equation} \label{fs25ssv}
\frac{\mathbb{P}(H_{n_0,n}  \geq x)}{1-\Phi \left( x\right)}=1+o(1)\ \ \ \ \textrm{and} \ \ \ \ \frac{\mathbb{P}(H_{n_0,n}  \leq -x)}{1-\Phi \left( x\right)}=1+o(1)
\end{equation}
uniformly for $0\leq x=o (  \sqrt{\log n}   ).$ The upper $(\kappa_n/2)$th quantile of a standard normal distribution satisfies
$$\Phi^{-1}( 1-\kappa_n/2) = O(\sqrt{| \ln \kappa_n |} ), $$  which, by (\ref{fs25sv}), is of order $ o (  \sqrt{\log n}   ).$
Proposition \ref{pro3.3} follows from  applying (\ref{fs25ssv})   to $H_{n_0,n}$. \hfill\qed

\subsection{An infectious disease model}
 An infectious disease model $(Z_n)_{n\geq 0}$ may be described as follows:
 \begin{equation}\label{ddfsfdgs}
 Z_0=1,\ \ \ \ Z_{n+1}= Z_{n}+\sum_{i=1}^{Z_n} Y_{n,i}, \   \ \ \ \textrm{for } n \geq 0,
 \end{equation}
where $Z_{n}$ stands for the total population of patients with infectious disease  at time $n$,
and $Y_{n,i}$ is the number of patients infected by the $i$-th individual of $Z_n$ in a unit time (for instance, one day). 
Moreover, we assume that  the random variables $  (Y_{n,i})_{i\geq 1} $ are iid random variables  with common  distribution law
  \begin{equation} \label{dssdfs}
 \mathbb{P}(Y_{n,i} =k     )     = p_k,\ \ \ k\in \mathbb{N},
\end{equation}
and are also independent to $Z_n.$   Denote by $r$  the average number of patients infected by an individual patient in a unite time,
  that is
  $$r=\mathbb{E} Y_{n,i}=\sum_{k=0}^ \infty k \, p_k .$$
   Denote by $v$ the standard variance of  $Y_{n,i}, n, i \geq 1$, then $v$ is also the standard variance of
  $ Z_1,$ that is   $$ v^2= \mathbb{E}(Z_1-m)^2.$$    To avid triviality, assume that $v>0.$ We are interested
  in the estimation of $r.$
\begin{proposition} \label{fdfgs}
Assume that $ \mathbb{E} Z_1 ^{2+\rho}< \infty$ for some $  \rho \in (0,  1]$. Let $\kappa_n \in (0, 1).$  Assume that
\begin{eqnarray}
 \big| \ln \kappa_n \big| =o \big( n^{\rho/(2+\rho)}  \big) .
\end{eqnarray}
Let $A_{n_0, n}$ and $B_{n_0, n}$ be defined in Proposition \ref{c0klfdgs}.
Then $[A_{n_0, n} -1 ,\     B_{n_0, n}-1 ]$ is a  $1-\kappa_n$ confidence interval for $r$, for $n$ large enough.
\end{proposition}
\emph{Proof.}  It is easy to see that (\ref{ddfsfdgs}) can be rewritten in the  form of (\ref{GWP}), with $X_{n,i}=1+Y_{n,i}$.
Thus, we have $m=1+r.$ Then Proposition \ref{fdfgs} follows by Proposition \ref{c0klfdgs}. \hfill\qed

 \section{Proof of Theorem \ref{th01}}
 \setcounter{equation}{0}
 In the proof of Theorem \ref{th01}, we will make use  of the following lemma (cf.\,Corollary 2.3 of Fan et al.\,\cite{FGLS19}), which
gives  self-normalized Cram\'{e}r moderate deviations for martingales.
 \begin{lemma}\label{lema01}
Let $(\eta_k, \mathcal{F}_k)_{k=1,...,n}$ be a finite sequence of martingale differences.
Assume that   there exist a constant $\rho \in (0, 1]$ and numbers $\gamma_n>0$ and $\delta_n\geq 0$ satisfying $\gamma_n, \delta_n \rightarrow 0$  such that
for all $1\leq i\leq n,$
\begin{equation} \label{cond3}
 \mathbb{E}[ |\eta_{k}|^{2+\rho} |\mathcal{F}_{k-1 } ]  \leq   \gamma_n^\rho \mathbb{E}[  \eta_{k} ^2 |\mathcal{F}_{k-1 } ]
\end{equation}
and
\begin{equation}\label{cond2}
 \Big\| \sum_{k=1}^n \mathbb{E}[  \eta_{k} ^2 |\mathcal{F}_{k-1 } ]-1\Big\|_\infty   \leq \delta_n^2  \ \  \ \ \textrm{a.s.}
\end{equation}
Denote
$$V_n= \frac{\sum_{k=1}^n   \eta_{k} }{ \sqrt{\sum_{k=1}^n   \eta_{k}^2} \ }$$
and
$$ \widehat{\gamma}_n(x, \rho) =  \frac{ \gamma_n^{ \rho(2-\rho)/4 } }{ 1+ x  ^{  \rho(2+\rho)/4 }}. $$
\begin{description}
  \item[\textbf{[i]}]  If $\rho \in (0, 1)$, then for all $0\leq x = o(  \gamma_n^{-1}  ),$
\begin{equation}
\bigg|\ln \frac{\mathbb{P}(V_n \geq x)}{1-\Phi \left( x\right)} \bigg| \leq      C_{\rho} \bigg( x^{2+\rho}  \gamma_n^\rho+ x^2 \delta_n^2 +(1+x)\Big(  \delta_n + \widehat{\gamma}_n(x, \rho) \Big) \bigg)   .
\end{equation}
  \item[\textbf{[ii]}]  If $\rho =1$, then for all $0\leq x = o(  \gamma_n^{-1}  ),$
\begin{equation}
\bigg|\ln \frac{\mathbb{P}(V_n \geq x)}{1-\Phi \left( x\right)} \bigg| \leq  C  \bigg( x^{3}  \gamma_n + x^2 \delta_n^2+(1+x)\Big(   \delta_n+ \gamma_n |\ln \gamma_n|  +\widehat{\gamma}_n(x, 1) \Big) \bigg)  .
\end{equation}
\end{description}
\end{lemma}

Now, we are in position to prove
Theorem \ref{th01}.
Denote $$\hat{\xi}_{k+1}= \sqrt{Z_{ k}} ( Z_{ k+1}/Z_{ k} -m ),$$
$\mathfrak{F}_{n_0} =\{ \emptyset, \Omega \}  $ and $\mathfrak{F}_{k+1}=\sigma \{ Z_{i}: n_0\leq i\leq k+1  \}$ for all $k\geq n_0$.
Notice that $X_{k,i}$ is independent of $Z_k.$
Then it is easy to verify that
 \begin{eqnarray}
\mathbb{E}[ \hat{\xi}_{k+1}  |\mathfrak{F}_{k } ] &=&  Z_{ k} ^{-1/2}  \mathbb{E}[   Z_{ k+1}  -mZ_{ k}  |\mathfrak{F}_{k } ] = Z_{ k} ^{-1/2} \sum_{i=1}^{Z_k} \mathbb{E}[   X_{ k, i}  -m   |\mathfrak{F}_{k } ]  \nonumber \\
& =&   Z_{ k} ^{-1/2} \sum_{i=1}^{Z_k} \mathbb{E}[   X_{ k, i}  -m   ] \nonumber \\
& =&0.
\end{eqnarray}
 Thus
 $(\hat{\xi}_k, \mathfrak{F}_k)_{k=n_0+1,...,n_0+n}$ is a finite sequence of martingale  differences.
Notice that   $X_{ k, i}-m, i\geq 1,$  are  centered and independent random variables. Thus,
 the following   equalities hold
\begin{eqnarray}
 \sum_{k=n_0}^{n_0+n-1} \mathbb{E}[ \hat{\xi}_{k+1}^2  |\mathfrak{F}_{k } ] &=&
\sum_{k=n_0}^{n_0+n-1}  Z_{k}^{-1} \mathbb{E}[  ( Z_{ k+1} -mZ_{ k} )^2  |\mathfrak{F}_{k } ]
 = \sum_{k=n_0}^{n_0+n-1}  Z_{k}^{-1} \mathbb{E}[  ( \sum_{i=1}^{Z_k}    (X_{ k, i}  -m)  )^2  |\mathfrak{F}_{k } ] \nonumber \\
&=&\sum_{k=n_0}^{n_0+n-1}  Z_{k}^{-1} Z_k \mathbb{E}[    (X_{ k, i}  -m)^2  ] \nonumber \  \\
&=& n v^2.  \label{ineq9.2}
\end{eqnarray}
Moreover, it is easy to see that
\begin{eqnarray}
  \mathbb{E}[ |\hat{\xi}_{k+1}|^{2+\rho} |\mathfrak{F}_{k } ] &=&
  Z_{k}^{-1-\rho/2} \mathbb{E}[  | Z_{ k+1} -mZ_{ k} |^{2+\rho }  |\mathfrak{F}_{k } ]  \nonumber \\
& = &  Z_{k}^{-1-\rho/2}   \mathbb{E}[  | \sum_{i=1}^{Z_k}    (X_{ k, i}  -m)  |^{2+\rho }  |\mathfrak{F}_{k } ] . \label{ines2.8}
\end{eqnarray}
By Rosenthal's inequality, we have
\begin{eqnarray*}
 \mathbb{E}[  | \sum_{i=1}^{Z_k}    (X_{ k, i}  -m)  |^{2+\rho }  |\mathfrak{F}_{k } ]
&\leq& C'_\rho \bigg(\Big(  \sum_{i=1}^{Z_k}\mathbb{E}( X_{ k, i}  -m )^2  \Big)^{1+\rho/2} + \sum_{i=1}^{Z_k}\mathbb{E}| X_{ k, i}  -m   |^{2+\rho } \bigg)  \\
& \leq& C'_\rho  \bigg(    Z_k^{1+\rho/2} v^{2+\rho} + Z_k  \mathbb{E}| Z_1  -m   |^{2+\rho } \bigg).
\end{eqnarray*}
Since the set of extinction of the process
$(Z_{k})_{k\geq 0}$ is negligible with respect to the annealed law $\mathbb{P}$, we have $Z_k\geq 1$ for any $k$.
 From (\ref{ines2.8}), by the last inequality  and the fact $Z_k\geq 1$, we deduce that
\begin{eqnarray}
  \mathbb{E}[ |\hat{\xi}_{k+1}|^{2+\rho} |\mathfrak{F}_{k } ] &\leq& C'_\rho  (     v^{ \rho} +   \mathbb{E}| Z_1  -m   |^{2+\rho }/v^2  )  v^2 \nonumber \\
& = & C'_\rho  (     v^{ \rho} +   \mathbb{E}| Z_1  -m   |^{2+\rho }/v^2  )    \mathbb{E}[ \hat{\xi}_{k+1}^{2} |\mathfrak{F}_{k } ]\nonumber \\
& = & C_\rho  (     v^{ \rho} +   \mathbb{E}  Z_1  ^{2+\rho }/v^2  )    \mathbb{E}[ \hat{\xi}_{k+1}^{2} |\mathfrak{F}_{k } ]. \label{ineq9.4}
\end{eqnarray}
Let $\eta_k =\hat{\xi}_{n_0+k}/\sqrt{n} v$ and $\mathcal{F}_{k}=\mathfrak{F}_{n_0+k}$. Then $(\eta_k,  \mathcal{F}_{k})_{k=1,...,n}$ is a  martingale difference sequences and
 satisfies the conditions  (\ref{cond3}) and (\ref{cond2}) with
$\delta_n=0$ and $\gamma_n=( C_\rho  (     v^{ \rho} +   \mathbb{E}  Z_1  ^{2+\rho }/v^2  )  )^{1/\rho}/ \sqrt{n} v$.
Clearly, it holds
$$  M_{n_0,n}=  \frac{ \sum_{k=1}^{n}  \eta_{k } }{\sqrt{\sum_{k=1}^{n} \eta_{k }^2 } \ }. $$
Applying Lemma \ref{lema01}  to $(\eta_k, \mathcal{F}_{k})_{k=1,...,n}$,   we obtain the desired inequalities. \hfill\qed

\section{Proof of Corollary \ref{coro01}}
\setcounter{equation}{0}
We first show that for any Borel set $B\subset \mathbb{R},$
\begin{eqnarray}\label{dfgsfdf}
\limsup_{n\rightarrow \infty}\frac{1}{a_n^2}\ln \mathbb{P}\bigg(\frac{M_{n_0,n} \ }{ a_n  }  \in B \bigg)   \leq  - \inf_{x \in \overline{B}}\frac{x^2}{2}.
\end{eqnarray}
When $B  =\emptyset,$ the last inequality is obvious, with   $-\inf_{x \in \emptyset}\frac{x^2}{2}=-\infty$. Thus,  we may assume that $B  \neq \emptyset.$ Let $x_0=\inf_{x\in B} |x|.$ Clearly, we have $x_0\geq\inf_{x\in \overline{B}} |x|.$
Then, by   Theorem \ref{th01}, it follows that for $a_n =o(\sqrt{n}),$
\begin{eqnarray*}
 \mathbb{P}\bigg(\frac{ M_{n_0,n} \ }{ a_n  }  \in B \bigg)
 &\leq&  \mathbb{P}\bigg( | M_{n_0,n}|  \geq a_n x_0\bigg)\\
 &\leq&  2\Big( 1-\Phi \left( a_nx_0\right)\Big)
  \exp\bigg\{ C_\rho    \bigg(    \frac{ ( a_nx_0)^{2+\rho} }{n^{\rho/2}}  + \frac{\ln n }{\sqrt{n} }+   \frac{ (1+a_nx_0)^{1-\rho(2+\rho)/4} }{n^{\rho(2-\rho)/8}}     \bigg)  \bigg\}.
\end{eqnarray*}
Using the following inequalities
\begin{eqnarray}\label{fgsgj1}
\frac{1}{\sqrt{2 \pi}(1+x)} e^{-x^2/2} \leq 1-\Phi ( x ) \leq \frac{1}{\sqrt{ \pi}(1+x)} e^{-x^2/2}, \ \   x\geq 0,
\end{eqnarray}
  and  the fact that $a_n \rightarrow \infty$ and $a_n/\sqrt{n}\rightarrow 0$,
we obtain
\begin{eqnarray*}
\limsup_{n\rightarrow \infty}\frac{1}{a_n^2}\ln \mathbb{P}\bigg(\frac{ M_{n_0,n}   \ }{ a_n  }  \in B \bigg)
 \ \leq \  -\frac{x_0^2}{2} \ \leq \  - \inf_{x \in \overline{B}}\frac{x^2}{2} ,
\end{eqnarray*}
which gives (\ref{dfgsfdf}).

Next, we prove that
\begin{eqnarray}\label{dfsffsfn}
\liminf_{n\rightarrow \infty}\frac{1}{a_n^2}\ln \mathbb{P}\bigg(\frac{ M_{n_0,n} \ }{ a_n   }  \in B \bigg) \geq   - \inf_{x \in B^o}\frac{x^2}{2} .
\end{eqnarray}
When $B^o =\emptyset,$ the last inequality is obvious, with   $ -\inf_{x \in  \emptyset}\frac{x^2}{2}=-\infty$. Thus, we may assume that $B^o \neq \emptyset.$ Since $B^o$ is an open set,
for any given small $\varepsilon_1>0,$ there exists an $x_0 \in B^o,$ such that
\begin{eqnarray*}
 0< \frac{x_0^2}{2} \leq   \inf_{x \in B^o}\frac{x^2}{2} +\varepsilon_1.
\end{eqnarray*}
Again by the fact that $B^o$ is an open set, for $x_0 \in B^o$ and all small enough $\varepsilon_2 \in (0, |x_0|], $ it holds $(x_0-\varepsilon_2, x_0+\varepsilon_2]  \subset B^o.$  Without loss of generality, we may assume that $x_0>0.$
Clearly, we have
\begin{eqnarray}
\mathbb{P}\bigg(\frac{M_{n_0,n} \ }{ a_n  } \in B  \bigg)   &\geq&   \mathbb{P}\bigg( M_{n_0,n}    \in (a_n ( x_0-\varepsilon_2), a_n( x_0+\varepsilon_2)] \bigg) \nonumber \\
&=&   \mathbb{P}\bigg( M_{n_0,n} \geq  a_n ( x_0-\varepsilon_2)   \bigg)-\mathbb{P}\bigg( M_{n_0,n}  \geq  a_n( x_0+\varepsilon_2) \bigg). \label{fsfvhd}
\end{eqnarray}
Again by Theorem  \ref{th01}, it is easy to see that for $a_n \rightarrow \infty$ and $ a_n =o(\sqrt{n} ),$ $$\lim_{n\rightarrow \infty} \frac{\mathbb{P}\big(M_{n_0,n} \geq  a_n( x_0+\varepsilon_2) \big) }{\mathbb{P}\big(  M_{n_0,n} \geq  a_n ( x_0-\varepsilon_2)   \big)  } =0 .$$
From (\ref{fsfvhd}),  by the last line and Theorem \ref{th01}, it holds for all $n$ large enough and $a_n =o(\sqrt{n} ),$
\begin{eqnarray*}
&&\mathbb{P}\bigg(\frac{M_{n_0,n} \ }{ a_n  } \in B  \bigg) \ \geq\    \frac12 \mathbb{P}\bigg(   M_{n_0,n}   \geq  a_n ( x_0-\varepsilon_2)   \bigg) \\
&& \ \ \ \ \ \ \ \ \  \ \geq   \   \frac12 \Big( 1-\Phi \left( a_n( x_0-\varepsilon_2)\right)\Big)   \exp\bigg\{ -C_\rho \bigg(    \frac{ ( a_nx_0)^{2+\rho} }{n^{\rho/2}} + \frac{\ln n }{\sqrt{n} } +   \frac{ (1+a_nx_0)^{1-\rho(2+\rho)/4} }{n^{\rho(2-\rho)/8}}     \bigg) \bigg\}.
\end{eqnarray*}
Using   (\ref{fgsgj1}) and  the fact that $a_n \rightarrow \infty$ and $a_n/\sqrt{n}\rightarrow 0$, after some calculations,
we get
\begin{eqnarray*}
 \liminf_{n\rightarrow \infty}\frac{1}{a_n^2}\ln \mathbb{P}\bigg(\frac{M_{n_0,n}   \ }{ a_n  }  \in B \bigg)  \geq  -  \frac{1}{2}( x_0-\varepsilon_2)^2 . \label{ffhms}
\end{eqnarray*}
Letting $\varepsilon_2\rightarrow 0,$  we  deduce that
\begin{eqnarray*}
\liminf_{n\rightarrow \infty}\frac{1}{a_n^2}\ln \mathbb{P}\bigg(\frac{ M_{n_0,n}  \ }{ a_n  }  \in B \bigg) \ \geq\ -  \frac{x_0^2}{2}  \  \geq \   -\inf_{x \in B^o}\frac{x^2}{2} -\varepsilon_1.
\end{eqnarray*}
Since that $\varepsilon_1$ can be arbitrarily small, we get (\ref{dfsffsfn}).
Combining (\ref{dfgsfdf}) and (\ref{dfsffsfn}) together, we complete  the proof of Corollary  \ref{coro01}.
 \hfill\qed

\section{Proof of Theorem \ref{jsdf}}
\setcounter{equation}{0}

In the proof of Theorem \ref{jsdf}, we will make use  of the following lemma (cf.\ Theorem 2.1 of Fan \cite{F19}), which
gives    exact  Berry-Esseen's bounds  for martingales.
\begin{lemma}\label{lema1}
Assume the conditions of Lemma \ref{lema01}.
\begin{description}
  \item[\textbf{[i]}]  If $\rho \in (0, 1)$, then
\begin{equation}
\sup_{ x \in \mathbb{R}}\Big|\mathbb{P}(\sum_{k=1}^n   \eta_{k}    \leq x)-\Phi \left( x\right) \Big|\leq C_{\rho} \Big( \gamma_n^\rho   + \delta_n \Big) .
\end{equation}
  \item[\textbf{[ii]}]  If $\rho =1$, then
\begin{equation}
\sup_{ x \in \mathbb{R}}\Big|\mathbb{P}( \sum_{k=1}^n    \eta_{k} \leq x)-\Phi \left( x\right) \Big| \leq C \, \Big( \gamma_n  |\log \gamma_n|  + \delta_n \Big).
\end{equation}
\end{description}
\end{lemma}

Recall the martingale  differences $(\eta_k, \mathcal{F}_{k})_{k=1,...,n }$  defined in the proof of  Theorem \ref{th01}.
Then  $\eta_k$ satisfies  the conditions  (\ref{cond3}) and (\ref{cond2}) with
$\delta_n=0$  and  $ \gamma_n=( C_\rho  (     v^{ \rho} +   \mathbb{E}  Z_1  ^{2+\rho }/v^2  )  )^{1/\rho}/ \sqrt{n} v.$
Clearly, it holds $H_{n_0,n}= \sum_{k=1}^{n}   \eta_k.$
Applying Lemma \ref{lema1} to $(\eta_k, \mathcal{F}_{k})_{k=1,...,n}$,  we obtain  the desired  inequalities. \hfill\qed

\section{Proof of Theorem \ref{fsfs}}
\setcounter{equation}{0}

Define the generating function of $Z_{n}$ as $f_{n}(s)=\mathbb{E} s^{Z_{n}},$ $ |s|\leq 1.$
We have the following lemma, see Athreya \cite{A94}.
\begin{lemma}\label{lemma01}
If $p_1>0$ then
\begin{equation}
\lim\limits_{n\to\infty}\frac{f_n(s)}{p_1^n}=\sum_{k=1}^{\infty}q_{k}s^k,
\end{equation}
where $(q_{k}, k\geq 1)$   is defined via the generating function $Q(s)=\sum_{k=1}^{\infty}q_{k}s^k, 0\leq s<1,$ the unique solution of the functional equation
$$Q(f(s))=p_{1}Q(s), \quad \mbox{ where } f(s) =\sum_{j=1}^{\infty}p_{j}s^j, \ 0\leq s<1,$$
subject to
$$Q(0)=0, \qquad  Q(1)=\infty, \qquad Q(s)<\infty \mbox{ for } 0\leq s <1.$$
\end{lemma}

\begin{lemma}\label{lemma2}
It holds
\begin{eqnarray}\label{fgffdgdfg}
\mathbb{P}(Z_n  \leq n)&\leq&  C_1 \exp\{  -n c_0\}.
\end{eqnarray}
\end{lemma}
\textit{Proof.}
When $p_1>0,$ using Markov's inequality and Lemma \ref{lemma01},  we have for $s_0=\frac{1+p_1}{2} \in (0, 1),$
\begin{eqnarray}
\sum_{k=1}^{n-1}\mathbb{P}(Z_n=k) I_{k}(x) &\leq&\mathbb{P}(Z_n  \leq n) = \mathbb{P}(s_0^{Z_n}\geq s_0^n) \leq  s_0^{-n}f_n(s_0) \nonumber  \\
&\leq& C (\frac{p_1}{s_0})^nQ(s_0) \nonumber \\
&=& C_1 \exp\{  -n c_0\},
\end{eqnarray}
where $C_1=C  Q(s_0) $ and $c_0= \ln (s_0/p_1).$ Notice that
$s_0  \in (p_1, 1),$ thus $c_0>0.$ Recall that $p_0=0.$
When $p_1=0,$ we have $Z_n\geq 2^n,$ and (\ref{fgffdgdfg}) holds obviously for all $n$ large enought.

In the proof of Theorem \ref{fsfs}, we need  the following technical lemma of Jing, Shao and Wang \cite{JSW03},
 which  gives a self-normalized Cram\'{e}r moderate deviation result for  iid random variables.
\begin{lemma}\label{crself}
Let $(Y_{i})_{i\geq 1} $ be a sequence of iid and  centered random variables. Assume that $\mathbb{E}|Y_1|^{2+\rho}< \infty$ for some $ \rho \in (0,1].$ Let $S_{n}=\sum_{i=1}^{n} Y_{i}$ and $V_{n}^{2}=\sum_{i=1}^nY_{i}^2$.
Then
\begin{equation}\label{ineq02}
\bigg|\ln\frac{\mathbb{P}(S_{n}/ V_{n}\geq x)}{1-\Phi(x)} \bigg| \leq  C_\rho  \frac{1+x^{2+\rho}}{ n^{\rho/2}}
\end{equation}
uniformly for $0 \leq x =o(\sqrt{n} )$ as $n\rightarrow \infty$.
\end{lemma}

\subsection{Proof of the theorem}
Now, we are in a position to prove  Theorem \ref{fsfs}.
Recalling that $Z_{n}$ is the number of individuals of the BPRE in generation $n$, and $X_{n,i}, \ 1\leq i\leq Z_{n},$
is the number of the offspring of the $i$th individual in generation $n$.
Denote
\begin{equation}
V(n)^2=\sum_{i=1}^{Z_{n}}(X_{n,i}-m )^2,\qquad \bar{X}(n)=\frac{  Z_{n+1}}{Z_{n}},\ \quad  \bar{Y}_n=\frac{  Z_{n+1}}{n}.
\end{equation}
Then we have
\begin{eqnarray}
\sum_{i=1}^{Z_n}(X_{n,i}-\bar{X}(n))^2&=& \sum_{i=1}^{Z_n}\big((X_{n,i}-m )+(m -\bar{X}(n) \big)^2  \nonumber \\
&=&V(n)^2-Z_n(m-\bar{X}(n))^2.\label{ineq10}
\end{eqnarray}
By   (\ref{ineq10}), it is easy to see that $T_n$ can be rewritten as follows:
$$ T_n  =  \frac{\sum_{i=1}^{Z_n}(X_{n,i}-m)}{  \sqrt{V(n)^2-Z_n(m-\bar{X}(n))^2}}  .$$
Notice that $X_{n,i}, \ 1\leq i\leq Z_{n},$ have the same distribution as $Z_1$, and that $Z_n$ is independent of $\xi_n$.   By the total probability formula and the independence of $Z_{n}$ and $(X_{n,i})_{i\geq 1}$, we obtain, for all $x\geq 0,$
\begin{eqnarray}
\mathbb{P} \Big( T_n \geq   x\Big) &=&\mathbb{P} \Bigg( \sum_{i=1}^{Z_n}(X_{n,i}-m) \geq x   \sqrt{V(n)^2-Z_n(m-\bar{X}(n))^2}\Bigg) \nonumber \\
&=&\sum_{k=1}^{\infty}\mathbb{P} (Z_n=k)\mathbb{P} \Bigg( \sum_{i=1}^k(X_{n,i}-m) \geq x    \sqrt{V_k^2-k(m-\bar{Y}_k)^2}\Bigg)\nonumber \\
&= &\sum_{k=1}^{\infty}\mathbb{P}(Z_n=k)  \mathbb{P} \Bigg( \sum_{i=1}^k(X_{n,i}-m) \geq x    \sqrt{V_k^2-k(m-\bar{Y}_k)^2}\Bigg) \nonumber \\
&=:&\sum_{k=1}^{\infty}\mathbb{P}(Z_n=k) I_{k}(x). \label{fgdg30}
\end{eqnarray}
By Lemma \ref{lemma01},  we have
\begin{eqnarray} \label{fgfdg30}
\sum_{k=1}^{n-1}\mathbb{P}(Z_n=k) I_{k}(x)  \leq \mathbb{P}(Z_n  \leq n)
 \leq   C_1 \exp\{  -n c_0\},
\end{eqnarray}
For $k\geq n,$ the tail probability  $I_{k}(x)$ can be divided into two parts: for all $x\geq0,$
\begin{eqnarray}
I_k(x)&=&\mathbb{P} \bigg(\sum_{i=1}^k(X_{n,i}-m)\geq x \sqrt{V_k^2-k(m-\bar{Y}_k)^2}, \ k(m-\bar{Y}_k)^2< V_k^2(1+x^\rho)/ k^{\rho/2} \bigg)\nonumber\\
&&+\ \mathbb{P} \bigg(\sum_{i=1}^k(X_{n,i}-m)\geq x \sqrt{V_k^2-k(m-\bar{Y}_k)^2}, \ k(m-\bar{Y}_k)^2\geq  V_k^2(1+x^\rho)/ k^{\rho/2} \bigg)\nonumber\\
&\leq& \mathbb{P} \bigg(\sum_{i=1}^k(X_{n,i}-m)\geq x V_k\sqrt{ 1-(1+x^\rho)/ k^{\rho/2}  }\bigg)+\mathbb{P} \bigg(k(m-\bar{Y}_k)^2 \geq    V_k^2(1+x^\rho)/ k^{\rho/2} \bigg) \nonumber  \\
&=:& I_{k,1}(x) +I_{k,2}(x). \label{thn303}
\end{eqnarray}

We first give an estimation for   $I_{k,1}(x).$ Notice that $(X_{n,i}-m)_{i\geq1}$ are conditional independent with respect to $\xi_n$.
When $k\geq n,$ by self-normalized moderate  deviations for   centered random variables $(X_{n,i}-m)_{i\geq1}$ (cf.\ Lemma \ref{crself}), we have, for all $0 \leq x =o(\sqrt{n}),$
$$\Bigg|\ln\frac{I_{k,1}(x)}{1-\Phi\big(x\sqrt{ 1-(1+x^\rho)/ k^{\rho/2}}  \ \big)} \Bigg| \leq  C_2 \frac{1+x^{2+\rho}}{ k^{\rho/2}} \leq C_2  \frac{1+x^{2+\rho}}{ n^{\rho/2}}.$$
Using  (\ref{fgsgj1}),
we deduce  that,  for all $x\geq 0$ and $0\leq \varepsilon   \leq 1$,
\begin{eqnarray}
  \frac{1-\Phi \left( x \sqrt{1- \varepsilon} \right)}{1-\Phi \left( x\right) }& =& 1+ \frac{ \int_{x  \sqrt{1-  \varepsilon}  }^x \frac{1}{\sqrt{2\pi}}e^{-t^2/2}dt }{1-\Phi \left( x\right) }
   \leq     1+ \frac{\frac{1}{\sqrt{2\pi}} e^{-x^2(1- \varepsilon)/2}  x \varepsilon  }{ \frac{1}{\sqrt{2 \pi} (1+x)} e^{-x^2/2}  }  \nonumber \\
   &\leq &   1+  C (1+  x^2)    \varepsilon  e^{     x^2 \varepsilon /2 }   \nonumber \\
    &\leq&     \exp\Big\{ C  (1+  x^2 )  \varepsilon  \Big\}. \label{sfdsh}
\end{eqnarray}
Using the last inequality, we get, for all  $k\geq n $  and all $0 \leq x =o(\sqrt{n}),$
\begin{eqnarray}
I_{k,1}(x)&\leq & \Big( 1- \Phi(x\sqrt{ 1-(1+x^\rho)/ k^{\rho/2}}  \ ) \Big)\exp\Big\{   C_2  \frac{1+x^{2+\rho}}{ n^{\rho/2}} \Big\}  \nonumber \\
&\leq& \Big( 1- \Phi(x) \Big)\exp\Big\{   C_2  \frac{1+x^{2+\rho}}{ n^{\rho/2}}+ C (1+x^2)\frac{ 1+x^\rho }{ k^{\rho/2}} \Big\} \nonumber \\
&\leq& \Big( 1- \Phi(x) \Big)\exp\Big\{   C_3  \frac{1+x^{2+\rho}}{ n^{\rho/2}}  \Big\} , \label{sdfgdsh1}
\end{eqnarray}
which gives an estimation for  $I_{k,1}(x).$

Next we give an estimation for   $I_{k,2}(x).$
Notice that
$$k(m-\bar{Y}_k)^2=\frac{1}{k}\bigg(\sum_{i=1}^k(X_{n,i}-m)\bigg)^2.$$
Thus,   we have
\begin{eqnarray*}
I_{k,2}(x) &=&\mathbb{P} \bigg(\Big(\sum_{i=1}^k(X_{n,i}-m)\Big)^2\geq  k^{1-\rho/2} V_k^2(1+x^\rho) \bigg)\nonumber\\
&=&\mathbb{P} \bigg( \Big|\sum_{i=1}^k(X_{n,i}-m)\Big| \geq V_k\sqrt{k^{1-\rho/2} (1+x^\rho)}\bigg) .
\end{eqnarray*}
Applying (\ref{ineq02})  to the centered random variables $(\pm (X_{n,i}-m))_{i\geq1}$, we obtain, for all $k\geq n$ and all $0\leq x =o(\sqrt{n})$,
\begin{eqnarray*}
I_{k,2}(x) &\leq& 2 \Big(1-\Phi(\sqrt{k^{1-\rho/2} (1+x^\rho)}\ )\Big)\exp\bigg\{ C\frac{1+(\sqrt{k^{1-\rho/2} (1+x^\rho)}\,)^{2+\rho}}{ \sqrt{k}}   \bigg \}\\
&\leq& 2 \exp\bigg\{ -\frac{1}{4}  k^{1-\rho/2} (1+x^\rho)   \bigg \},
\end{eqnarray*}
where the last line follows by (\ref{fgsgj1}).
Again by (\ref{fgsgj1}), we have, for all $k\geq n$ and all $0\leq x =o(\sqrt{n})$,
\begin{eqnarray}
I_{k,2}(x)  &\leq& 2 \exp\bigg\{ -\frac{1}{4}  n^{1-\rho/2} (1+x^\rho)   \bigg \} \nonumber \\
 &\leq&  C\frac{1+x}{n} \Big(1-\Phi(x)\Big), \label{sdfgdsh2}
\end{eqnarray}
which gives an estimation for  $I_{k,2}(x).$

Combining (\ref{thn303}), (\ref{sdfgdsh1}) and (\ref{sdfgdsh2}) together, we get, for all $k\geq n$ and all $0\leq x =o(\sqrt{n})$,
\begin{eqnarray}
I_k(x)&\leq& \Big( 1- \Phi(x) \Big)\exp\Big\{   C_4  \frac{1+x^{2+\rho}}{ n^{\rho/2}}  \Big\}.
\end{eqnarray}
Returning to (\ref{fgdg30}), using the last inequality and (\ref{fgfdg30}), we  deduce that, for all  $0\leq x =o(\sqrt{n} )$,
\begin{eqnarray}
\mathbb{P} \Big(T_n \geq   x\Big) &\leq& \sum_{k=1}^{n-1}\mathbb{P}(Z_n=k)I_{k}(x)   + \sum_{k=n}^{\infty}\mathbb{P}(Z_n=k) I_{k}(x)\nonumber\\
&\leq&  C_1 \exp\{  -C_0 n  \}   + \sum_{k=n}^{\infty}\mathbb{P}(Z_n=k) \Big( 1- \Phi(x) \Big)\exp\Big\{   C_4  \frac{1+x^{2+\rho}}{ n^{\rho/2}}  \Big\}\nonumber\\
&\leq&  C_1 \exp\{  -C_0 n  \} + \sum_{k=1}^{\infty}\mathbb{P}(Z_n=k) \Big( 1- \Phi(x) \Big)\exp\Big\{   C_4  \frac{1+x^{2+\rho}}{ n^{\rho/2}}  \Big\} \nonumber \\
&=&  C_1 \exp\{  -C_0 n \}  +   \Big( 1- \Phi(x) \Big)\exp\Big\{   C_4 \frac{1+x^{2+\rho}}{ n^{\rho/2}}  \Big\} \nonumber \\
&\leq&    \Big( 1- \Phi(x) \Big)\exp\Big\{   C_5  \frac{1+x^{2+\rho}}{ n^{\rho/2}}  \Big\}, \label{ineq002}
\end{eqnarray}
where the last line follows by (\ref{fgsgj1}).

Next, we consider the lower bound of $\mathbb{P}\big(T_n \geq   x\big).$
For $I_{k}(x)$, we have the following estimation: for all $k\geq n$ and  all $0\leq x =o(\sqrt{n}),$
\begin{eqnarray}
I_k(x)&=&\mathbb{P} \bigg( \sum_{i=1}^k(X_{n,i}-m) \geq x    \sqrt{V_k^2-k(m-\bar{Y}_k)^2}\bigg) \nonumber \\
&\geq& \mathbb{P} \bigg(\sum_{i=1}^k(X_{n,i}-m)\geq x V_k   \bigg).
\end{eqnarray}
When $k\geq n,$ by self-normalized moderate  deviations for iid random variables (cf.\ Lemma \ref{crself}), we have, for all $0 \leq x =o(\sqrt{n}),$
$$I_{k }(x) \geq \Big(1- \Phi(x) \Big) \exp\Big\{-C_6  \frac{1+x^{2+\rho}}{ n^{\rho/2}} \Big\}.$$
Returning to (\ref{fgdg30}),  we  deduce that, for all  $0\leq x =o(\sqrt{n})$,
\begin{eqnarray*}
\mathbb{P} \Big(T_n \geq   x\Big) &\geq&  \sum_{k=n}^{\infty}\mathbb{P}(Z_n=k) I_{k}(x)\\
&\geq& \Big(1- \Phi(x) \Big) \exp\Big\{-C_6  \frac{1+x^{2+\rho}}{ n^{\rho/2}} \Big\} \sum_{k=n}^{\infty}\mathbb{P}(Z_n=k)\\
&\geq& \Big(1- \Phi(x) \Big) \exp\Big\{-C_6 \frac{1+x^{2+\rho}}{ n^{\rho/2}} \Big\} \Big( 1-  \mathbb{P}(Z_n \leq n)\Big) .
\end{eqnarray*}
Using   Lemma \ref{lemma2},  we get,  for all  $0\leq x =o(\sqrt{n})$,
\begin{eqnarray}
\mathbb{P} \Big(T_n \geq   x\Big)
&\geq& \Big(1- \Phi(x) \Big) \exp\Big\{-C_6  \frac{1+x^{2+\rho}}{ n^{\rho/2}} \Big\} \Big( 1- C_1 e^{-C_0 n  }\Big) \nonumber \\
&\geq& \Big(1- \Phi(x) \Big) \exp\Big\{-C_9  \frac{1+x^{2+\rho}}{ n^{\rho/2}} \Big\}.  \label{ineq001}
\end{eqnarray}
Combining (\ref{ineq002}) and (\ref{ineq001}) together, we obtain the desired inequality.

Applying  (\ref{fsdcv01}) to $\big(m-X_{n,k}\big)_{ k \geq 1}$, we find that
(\ref{fsdcv01}) remains valid when $\frac{\mathbb{P} (T_n \geq x)}{1-\Phi \left( x\right)}$ is replaced  by $\frac{\mathbb{P} (T_n \leq -x)}{ \Phi \left( -x\right)}$. This completes the proof of Theorem \ref{fsfs}. \hfill\qed

\section{Proof of Corollary \ref{corollary03}}
\setcounter{equation}{0}
Clearly, it holds
\begin{eqnarray}
&& \sup_{ x \in \mathbb{R}   }  \Big|\mathbb{P} \big(  T_n \leq x  \big) -  \Phi \left( x\right) \Big| \nonumber  \\
 && \leq  \sup_{   x > n^{\rho/(8+4\rho)}  } \Big|\mathbb{P} \big(  T_n \leq x  \big) -  \Phi \left( x\right) \Big|    + \sup_{  0 \leq x \leq  n^{\rho/(8+4\rho)}  } \Big|\mathbb{P} \big(  T_n \leq x  \big) -  \Phi \left( x\right) \Big| \nonumber\\
&& \ \ \ \ \ +  \sup_{ -  n^{\rho/(8+4\rho)}\leq x \leq 0 } \Big|\mathbb{P} \big(  T_n \leq x  \big) -  \Phi \left( x\right) \Big|  +\sup_{   x < -  n^{\rho/(8+4\rho)}} \Big|\mathbb{P} \big(  T_n \leq x  \big) -  \Phi \left( x\right) \Big| \nonumber  \\
&&=:  TH_1 + TH_2+TH_3+TH_4.   \label{ineq0d10}
\end{eqnarray}
By Theorem \ref{fsfs} and (\ref{fgsgj1}), it is easy to see that
\begin{eqnarray*}
TH_1 &= & \sup_{   x >  n^{\rho/(8+4\rho)}} \Big|\mathbb{P} \big(  T_n> x  \big)- \big(1 -  \Phi \left( x\right) \big) \Big| \\
 &\leq & \sup_{   x >  n^{\rho/(8+4\rho)}}  \mathbb{P} \big(  T_n > x  \big) + \sup_{   x >  n^{\rho/(8+4\rho)}}  \big(1 -  \Phi \left( x\right) \big)   \\
 &\leq &   \mathbb{P} \big( T_n >  n^{\rho/(8+4\rho)}   \big) +    \big(1 -  \Phi (  n^{\rho/(8+4\rho)} ) \big) \\
 &\leq &  \big(1 -  \Phi (  n^{\rho/(8+4\rho)} ) \big)e^C   +    \exp\Big\{ -\frac12  (  n^{\rho/(8+4\rho)} )^2 \Big\} \\
 & \leq & \frac{ C_1 }{ n^{\rho/2}}
\end{eqnarray*}
and
\begin{eqnarray*}
TH_4  &\leq & \sup_{   x < -  n^{\rho/(8+4\rho)}}  \mathbb{P} \big(  T_n \leq x  \big) + \sup_{   x < -  n^{\rho/(8+4\rho)}}    \Phi \left( x\right)     \\
 &\leq &   \mathbb{P} \big( T_n \leq  -n^{\rho/(8+4\rho)}   \big) +     \Phi ( - n^{\rho/(8+4\rho)} )   \\
 &\leq &   \Phi ( - n^{\rho/(8+4\rho)} )  e^C   +    \exp\Big\{ -\frac12  (  n^{\rho/(8+4\rho)} )^2 \Big\} \\
 & \leq & \frac{ C_2 }{ n^{\rho/2}}.
\end{eqnarray*}
By Theorem \ref{fsfs} and the inequality $|e^x-1|\leq |x|e^{|x|},$ we have
\begin{eqnarray*}
TH_2 &=&\sup_{0\leq   x \leq n^{\rho/(8+4\rho)}} \Big|\mathbb{P} \big(  T_n > x  \big)- \big(1 -  \Phi \left( x\right) \big) \Big| \nonumber \\
  &\leq & \sup_{0\leq   x \leq n^{\rho/(8+4\rho)}} \Big(1-\Phi(x) \Big) \Big| e^{   C  ( 1+ x^{2+\rho})  / n^{\rho/2} }   -1 \Big| \nonumber\\
  &\leq & \frac{C}{n^{\rho/2}} \sup_{0\leq   x \leq n^{\rho/(8+4\rho)}}    \Big(1-\Phi(x) \Big) ( 1+ x^{2+\rho}) e^{   C  ( 1+ x^{2+\rho})  / n^{\rho/2} }  \nonumber\\
 &\leq &  \frac{C_3}{n^{\rho/2}}
 \end{eqnarray*}
 and, similarly,
 \begin{eqnarray*}
TH_3 &=& \sup_{ - n^{1/8}\leq x \leq 0 } \Big|\mathbb{P} \big(   T_n \leq x  \big) -  \Phi \left( x\right) \Big| \nonumber \\
  &\leq & \sup_{ - n^{1/8}\leq x \leq 0 } \Phi(x) \Big| e^{   C  ( 1+ |x|^3) (\ln n)/ \sqrt{n} }   -1 \Big| \nonumber\\
 &\leq & \frac{C_4}{n^{\rho/2}}.
 \end{eqnarray*}
Applying  the bounds  of $TH_1, TH_2, TH_3$ and $TH_4$ to (\ref{ineq0d10}),   we obtain the desired   inequality. This completes the proof of Corollary \ref{corollary03}. \hfill\qed

\section*{Acknowledgements}

This work has been partially supported by the National Natural Science Foundation
of China (Grant Nos.\,11601375 and 11971063).

\end{document}